\documentstyle[syntonly]{amsart}

\numberwithin{equation}{section}

\newtheorem{thm}{Theorem}[section]
\newtheorem{lem}{Lemma}[section]
\newtheorem{cor}{Corollary}[section]
\newtheorem{prop}{Proposition}[section]

\newtheorem{PBWthm}{PBW-Theorem}

\theoremstyle{definition}

\newtheorem{defn}{Definition}[section]
\newtheorem{exmp}{Example}[section]

\theoremstyle{remark}

\newtheorem{rem}{Remark}

\newcommand{\su}{U_q{\frak su}(1,1)}
\newcommand{\g}{U_q{\frak g}}
\newcommand{\h}{{\cal H}_{\frak g}(V)}
\newcommand{\F}{{\Bbb C}[V_{\Bbb R}]_q}
\newcommand{\C}{{\Bbb C}}
\newcommand{\e}{\varepsilon}
\newcommand{\f}{\varphi}
\newcommand{\T}{{\cal T}(V\oplus V^*)}
\newcommand{\usln}{U_q{\frak sl}(n+1)}
\newcommand{\usl}{U_q{\frak sl}(2)}
\newcommand{\R}{\check{R}}
\newcommand{\z}{\hat{z}}
\newcommand{\y}{\hat{y}}

\newcommand{\inv}{{\Bbb C}[{\cal O}_V]_q^{\rm inv}}
\newcommand{\Ch}{{\Bbb C}[{\cal O}_V]_q^{\frak h}}
\newcommand{\D}{\Delta}
\newcommand{\ra}{\rightarrow}
\newcommand{\mra}{\mapsto}
\newcommand{\om}{\omega}
\newcommand{\gc}{U_q{\frak g}_c}
\newcommand{\go}{U_q{\frak g}_0}
\newcommand{\ri}{@>\sim>>}
\newcommand{\hc}{{\cal H}_{{\frak g}_c}(V)}
\newcommand{\ho}{{\cal H}_{{\frak g}_0}(V)}
\newcommand{\sun}{U_q{\frak su}(n+1)}
\newcommand{\usu}{U_q{\frak su}\langle{\boldsymbol \iota}\rangle}
\newcommand{\CO}{{\Bbb C}[{\cal O}_V]_q}
\newcommand{\FX}{Func(X)_q}
\newcommand{\Fun}{Func({\cal O}_V)_q}
\newcommand{\cd}{Func(X_{c_0,d_0})_q}
\newcommand{\x}{\boldsymbol x}

\begin{document}

\title[On Quantum Orbit Method]
{Complimentary Series Representations And
Quantum Orbit Method}

\author{Leonid I. Korogodsky}

\subjclass{Quantum groups}


\address{Institute for Advanced Study, Princeton, NJ 08540, U.S.A.}

\email{leokor@@math.ias.edu}

\begin{abstract}

A version of quantum orbit method is presented for real forms of
equal rank of quantum complex simple groups. A quantum moment map
is constructed, based on the canonical isomorphism between a quantum
Heisenberg algebra and an algebra of functions on a family of quantum
$G$-spaces. For the series $A$, we construct some irreducible 
$*$-representations of $U_q{\frak g}$ which correspond to the semi-simple
dressing orbits of minimal dimension in the dual Poisson Lie group.
It is shown that some complimentary series representations correspond 
to some quantum 'tunnel' $G$-spaces which do not have a quasi-classical 
analog.

\end{abstract}

\maketitle

\section{Introduction}

This articles lies on the intersection of two areas of mathematics.
One is based on the tradition of the orbit method, pioneered by
A. Kirilov and B. Kostant, to realize the unitary irreducible
representations of Lie groups geometrically in certain bundles
on the orbits of the coadjoint action. And the other one is the
legacy of the V. Drinfeld's approach to quantum groups, which 
teaches us to look at the Lie groups as quasi-classical analogs
of their quantum counterparts.

An important observation that has given rise to various attempts
to develop a quantum analog of the classical orbit method is
that the coadjoint action is actually a special case of the 
dressing action of Poisson Lie groups. Namely, if $G$ is a
Poisson Lie group with the trivial Poisson structure, the
dual Poisson Lie group is isomorphic to the space $\g^*$, 
dual to its Lie algebra, which is considered as an Abelian
Lie group with the Poisson structure defined in terms of the
Kirillov-Kostant bracket.

Thus, one can look at the coadjoint action as a part of a more
general picture. In particular, if the Poisson Lie group structure
arises in the quasi-classical limit from a quantum group, we can
rise a question on a geometric realization of the irreducible
$*$-representations of the quantum universal enveloping algebra
in a way similar to how the classical orbit method works. It is
especially intriguing, for the classical orbit method doesn't
allow to realize all representations. For example, already in
the case of $SU(1,1)$ the complimentary series representations
cannot be realized in the classical case.

The quantum analog of the dressing action is based on the
fundamental principle of the Drinfeld's duality. Namely, the
same quantum algebra, say, the quantum universal enveloping
algebra, has two quasi-calssical limits. One being the classical
universal enveloping algebra $U_q\g$, and the other one being the
appropriate algebra of functions on the dual Poisson Lie group $G^*$.
A surprising manifestation of this phenomenon is that the
quantum analog of the dressing action becomes the quantum
adjoint action of our quantum algebra on itself, if we
assume that the copy of $U_q\g$ that serves as the
space of representation becomes the algebra of functions on $G^*$
in the limit, while the copy that acts on it becomes the
enveloping algebra instead.

As another manifestation of this duality, we construct a family
of quantum Heisenberg algebras which can be considered, at the same
time, as a family of quantum algebras of functions on the generalized
flag manifolds for $G$. The usual morphism from $U_q\g$ to such
quantum Heisenberg algebra, whose quasi-classical analog accounts
for the realization of the elements of $U\g$ in terms of the
differential operators, becomes, when considered as a morphism
from the quantum algebra of functions on $G^*$ into a quantum
algebra of functions on a $G$-space, a quantum analog of the 
corresponding moment map.

In this paper we proceed to introduce quantum analogs of polarizations
on the quantum $G$-spaces that arise from our construction. And eventually,
we provide a geometric realization of some series of the irreducible
$*$-representations on $U_q\g$. Namely, for the series $A$, the
representations corresponding to the minimal non-zero dimension orbits 
of the dressing action are realized geometrically.

It is worthy to note that in a hallmark example of the quantum group
$SU(1,1)$ all irreducible $*$-representations are realized that way,
including the complimentary and ``strange'' series representations,
the latter being specific for the quantum case. In particular, the
quantum $G$-spaces that account for the complimentary series
representations have a peculiar nature that we describe informally
as a 'tunnel' $G$-space. In the case of $SU(1,1)$ they are quantum
analogs of two-sheet hyperboloids, but they behave as if they were
connected. So that the quantum Haar measure (the invariant integral
on the quantum algebra of functions), when restricted to the subalgebra
of spherical functions, has a geometric progression for the support
that 'jumps' over the gap between the two 'sheets'.

\section{Quantum Heisenberg Algebras}

Throughout the paper we suppose that $q$ is real, $0<|q|<1$.

Suppose that $G$ is a finite-dimensional complex simple Poisson Lie 
group, $\frak g$ its Lie bialgebra, and $\g$ the corresponding 
quantized universal enveloping algebra (cf. \cite{D1,J1}). 
Let $V$ be a finite-dimensional simple $\g$-module, and $V^*$
the dual $\g$-module defined by
\[ \langle\xi\f,v\rangle=\langle\f,S(\xi)v\rangle, \]
where $\xi\in\g,~v\in V,~\f\in V^*$, and $S$ is the antipode
in $\g$. Then, we can define a {\em quantum Heisenberg algebra\/}
$\h$ as follows.

Let $R$ be the quantum $R$-matrix acting in $(V\oplus V^*)^{\otimes 2}$,
and $\R=PR$, where $P:(V\oplus V^*)^{\otimes 2}\ra(V\oplus V^*)
^{\otimes 2}$ is the usual permutation operator $a\otimes b\mra
b\otimes a$. As is well known, the operator $\R$ is invertible
and diagonalizable, has real spectrum and commutes with the
action of $\g$.

Consider the algebra $\F$ which is the quotient of the tensor
algebra 
\[ \T=\C\oplus(V\oplus V^*)\oplus
(V\oplus V^*)^{\otimes 2}\oplus\dotsb \]
over the two-sided ideal $J(W)$ generated by the span $W\subset
(V\oplus V^*)^{\otimes 2}$ of eigen-vectors of $\R$ with negative
eigen-values. 

The tensor algebra $\T$ has a canonical
$\g$-module algebra structure, which means that the canonical
$\g$-module structure defined by the action of $\g$ on
$V\oplus V^*$ is compatible with the algebra structure
so that the multiplication map $\T\otimes\T\ra\T$
is a morphism of $\g$-modules.

Since $\R$ commutes with the $\g$-action, the two-sided ideal
$J$ is a $\g$-submodule. This follows that $\F$ has a canonical
$\g$-module algebra structure as well. Note that $\F$ can be thought
of as a quantum analogue of the algebra of polynomial functions on
the $G$-space $V\oplus V^*$, the subalgebra generated
by $V^*$ (resp. $V$) playing the role of the algebra of holomorphic 
(resp. anti-holomorphic) polynomials. 

In the classical case this algebra has a canonical central 
extension given by
\[ ab-ba=\langle a,b \rangle C, \]
where $C$ is the central element, and $\langle \; , \; \rangle$
is the canonical bilinear form on $V\oplus V^*$ -- it is the
natural pairing between $V$ and $V^*$ and zero on both $V$ and $V^*$.
The result is known as the Heisenberg algebra. The quantum analogue 
is described below.

Consider the subspace $I$ of $\g$-invariant elements in $V\otimes V^*\oplus 
V^*\otimes V$ (that is, the elements $v$ such that $\xi v=\e(\xi)v$ 
for any $\xi\in\g$, where $\e$ is the counit). It is obviously 
two-dimensional, one generator in $V\otimes V^*$ and another one
in $V^*\otimes V$. Since $\R$ commutes with the action of
$\g$, $I$ is invariant with respect to $\R$. Since $\R$
permutes $V\otimes V^*$ and $V^*\otimes V$, it must have two
distinct eigen-values in $I$. From the other side, it is easy
to see that $\R^2|_I$ must be a constant, so that the
eigen-values must be of opposite sign. Thus, we have proved the
following proposition.

\begin{prop}
\label{one-dim}
The vector space $I_0=I\cap W$ is one-dimensional.
\end{prop}

Note that the $\g$-module $W$ is completely reducible, so that
there is a unique $\g$-submodule $W_0\subset W$ such that
$W=W_0\oplus I_0$. Let $J(W_0)\subset\T$
be the two-sided ideal generated by $W_0$.

\begin{defn}
\label{Heisenberg}
The $\g$-module algebra $\h=\T/J(W_0)$ is called {\em the
quantum Heisenberg algebra\/} corresponding to $\frak g$
and $V$. The following diagram is exact, where $C$ is the
image of a generator of $I_0$ and $p$ is the quotient map:
\[ 0\ra\C[C]\hookrightarrow \h @>p>> \F\ra 0. \]
\end{defn}

\begin{rem}
\label{Weyl}
The above definition first appeared in my preprint
\cite{K2}. Later we will observe how it is related to
the quantum Weyl algebra introduced in \cite{H}.
\end{rem}

\begin{exmp}
\label{exmp-1}
Suppose that ${\frak g}={\frak sl}(n+1)$ equipped with the
so-called {\em standard\/} Lie bialgebra structure. It is 
defined by the Manin triple $({\frak g}\oplus {\frak g},
{\frak g}, {\frak n}_+\oplus {\frak h}\oplus {\frak n}_-)$, 
where $\frak g$ is embedded into ${\frak g}\oplus {\frak g}$
as the diagonal, ${\frak n}_+$ as $({\frak n}_+,0)$,
${\frak n}_-$ as $(0,{\frak n}_-)$, $\frak h$ as
$\{(a,-a) \; | \; a\in {\frak h}\}$. Here $\frak h$ 
is the Cartan subalgebra of diagonal matrices, and
${\frak b}_+$ (resp. ${\frak b}_-$) the nilpotent subalgebra
of nilpotent upper- (resp. lower-) triagonal matrices.

Recall that $\usln$ is generated by $E_i,F_i,
K_i,K_i^{-1}$, $i=1,2,\dotsc,n$ with
the relations
\begin{gather}
\nonumber
\left[E_i,F_j\right]=\delta_{ij}\frac{K_i-K_i^{-1}}{q^{-1}-q}, 
 \; \; \; K_iK_j=K_jK_i, \\
\nonumber
K_iE_j=q^{-a_{ij}}E_jK_i, \; \; \; K_iF_j=q^{a_{ij}}F_jK_i, \\
\nonumber
E_iE_j=E_jE_i, (|i-j|>1), \; \; \; F_iF_j=F_jF_i, (|i-j|>1), \\
\nonumber
E_i^2E_{i+1}-\left(q+q^{-1}\right)E_iE_{i+1}E_i+E_{i+1}E_i^2=0, \\
\nonumber
F_i^2F_{i+1}-\left(q+q^{-1}\right)F_iF_{i+1}F_i+F_{i+1}F_i^2=0,
\end{gather}
where $a_{ii}=2$, $a_{i,i\pm 1}=-1$, and $a_{ij}=0$ otherwise.

The Hopf algebra structure on $\usln$ is given by
\begin{gather}
\nonumber
\D(K_i)=K_i\otimes K_i, \\
\nonumber
\D(E_i)=E_i\otimes 1+K_i^{-1}\otimes E_i, \; \; \;
\D(F_i)=F_i\otimes K_i+1\otimes F_i, \\
\nonumber
S(E_i)=-K_iE_i, \; \; \; S(F_i)=-F_iK_i^{-1},
\; \; \; S(K_i)=K_i^{-1}, \\
\nonumber
\e(E_i)=0, \; \; \; \e(F_i)=0, \; \; \; 
\e(K_i)=1,
\end{gather}
where $\D$ is the comultiplication, $S$ the antipode, and
$\e$ the counit.

Let $V$ be the finite-dimensional simple
$\usln$-module corresponding to the first fundamental
weight $\om_1$ (that is, the defining representation).
The algebra $\F$ in this case is
generated by $z_0,z_1,\dotsc,z_n$ and $\z_0,\z_1,\dotsc,
\z_n$ with the relations (cf. \cite{RTF})
\begin{eqnarray}
\nonumber
z_iz_j=qz_jz_i, \; (i<j), & \z_i\z_j=q^{-1}\z_j\z_i, \; (i<j), \\
\nonumber
z_i\z_j=q\z_jz_i, \; (i\not =j), & 
z_i\z_i-\z_iz_i=(q^{-2}-1)\sum_{k>i}z_k\z_k. 
\end{eqnarray}
Here $z_0,z_1,\dotsc,z_n$ and $\z_0,\z_1,\dotsc,
\z_n$ are nothing but the projections of the vectors
of the dual canonical bases of $V$ and $V^*$ respectively.

This means that the action of $\usln$ on $\F$ is given by
\begin{eqnarray}
\label{action1}
E_i:z_j\mra\delta_{ij}z_{j-1}, &
\z_j\mra -\delta_{i-1,j}q^{-1}\z_{j+1}, \\
\label{action2}
F_i:z_j\mra\delta_{i-1,j}z_{j+1}, &
\z_j\mra -\delta_{ij}q\z_{j-1}, \\
\label{action3}
K_i:z_j\mra
\begin{cases}
q^{-1}z_{i-1} & \text{if $j=i-1$,} \\
qz_i          & \text{if $j=i$,} \\
z_j           & \text{if otherwise,} 
\end{cases}, & \z_j\mra
\begin{cases}
q\z_{i-1}  & \text{if $j=i-1$,} \\
q^{-1}\z_i & \text{if $j=i$,} \\
\z_j       & \text{if otherwise.}
\end{cases}
\end{eqnarray}

In this case the subspace $I_0$ is spanned by
\[ \sum_{k=0}^nz_k\otimes\z_k-\sum_{k=0}^nq^{-2k}\z_k\otimes z_k. \]

The quantum Heisenberg algebra $\h$ is generated by
$z_0,z_1,\dotsc,z_n,\z_0,\z_1,\dotsc$, $\z_n$ and $C$
with the relations (cf. \cite{K2})
\begin{eqnarray}
\nonumber
z_iz_j=qz_jz_i, \; (i<j), & \z_i\z_j=q^{-1}\z_j\z_i, \; (i<j), \\
\nonumber
z_i\z_j=q\z_jz_i, \; (i\not =j), & 
z_i\z_i-\z_iz_i=C+\left(q^{-2}-1\right)\sum_{k>i}z_k\z_k, \\
\nonumber
z_i C=q^2C z_i, & \z_i C=q^{-2}C\z_i. 
\end{eqnarray}
The action of $\usln$ on $\h$ is given by \eqref{action1}-
\eqref{action3} and by 
\begin{equation}
\nonumber
\xi \, C=\e(\xi)C
\end{equation}
for any $\xi\in\usln$ (i.e., $C$ is a
$\usln$-invariant element).
\end{exmp}

Now suppose that ${\frak g}_c$ is a compact real form of $\frak g$,
which is unique up to an inner automorphism. Then, there is an 
antilinear involutive automorphism $\om$ of $\frak g$ such that
\begin{equation}
\nonumber
{\frak g}_c=\{a\in{\frak g}\; | \; \om(a)=a\}. 
\end{equation}

As is well known, there is a quantization of ${\frak g}_c$ in the
from of a Hopf $*$-algebra $\gc=(\g,\flat)$, where $\flat$ is
an antilinear involutive algebra anti-automorphism and coalgebra
automorphism such that
\begin{equation}
\nonumber
\om_c:\xi\mra (S(\xi))^{\flat}
\end{equation}
is an involution (thus, an antilinear involutive algebra automorphism 
and coalgebra anti-automorphism). Then, any finite-dimensional 
$\g$-module -- in particular, our module $V$ -- has a Hilbert space
structure which makes the action of $\gc$ into a $*$-representation.
Let $\iota:V \ri V^*$ be the antilinear isomorphism of vector spaces 
induced by the scalar product on $V$. The following proposition is
rather obvious.

\begin{prop}
\nonumber
The map $\iota:V\ri V^*$ can be uniquely extended to an antilinear 
involutive anti-automorphism $\sharp$ of $\h$ such that $C^{\sharp}=C$. 
Then, the $*$-algebra $\hc=(\h,\sharp)$ is a $\gc$-module $*$-algebra,
which means that for any $\xi\in\gc, f\in\hc$ one has 
\begin{equation}
\nonumber
(\xi f)^{\sharp}=\om_c(\xi)f^{\sharp}
\end{equation}
\end{prop}
(The definition of $\frak A$-module $*$-algebra -- where $\frak A$
is a Hopf $*$-algebra -- can be found in \cite{K1,K2,KV1}.)

Now suppose that $\go$ is a Hopf $*$-algebra $(\g,\natural)$
which is a quantization of a non-split real form ${\frak g}_0$
of $\frak g$ equipped with a standard Lie bialgebra structure
given by the Manin triple
\begin{equation}
\nonumber
({\frak g},{\frak g}_0,{\frak n}_+\oplus i{\frak h}_0),
\end{equation}
where ${\frak h}_0={\frak h}\cap {\frak g}_0$ (it depends,
of course, on the choice of the maximal nilpotent subalgebra
${\frak n}_+$). Let
\begin{equation}
\nonumber
\om_0:\xi\mra (S(\xi))^{\natural}
\end{equation}
be the corresponding antilinear involutive algebra automorphism
and coalgebra anti-automorphism on $\go$.

\begin{prop}
\nonumber
There exists a unique antilinear involutive anti-automorphism
$*$ of $\h$ such that $\ho=(\h,*)$ is a $\go$-module $*$-algebra,
which means that for any $\xi\in\go, f\in\ho$, one has 
\begin{equation}
\nonumber
(\xi f)^*=\om_0(\xi)f^*.
\end{equation}
\end{prop}

\begin{pf}
The composition $\tau=\om_0\om_c$ is a linear Hopf algebra
automorphism of $\g$. Then, there exists an operator $t:V
\ra V$ such that
\begin{equation}
\nonumber
(\tau\xi)(f)=(t^{-1}\xi t)(f)
\end{equation}
for any $\xi\in\g, f\in V$. The $\go$-module $*$-algebra 
structure on $\ho$ is given by
\begin{equation}
\nonumber
f^*=t(f^{\sharp}) \text{ for any } f\in V, \; \; C^*=C.
\end{equation}
The uniqueness is obvious.
\end{pf}

\begin{exmp}
\label{exmp-2}
Within the framework of Example \ref{exmp-1} we have that
the compact real form of $\usln$ is $\sun=(\usln,\flat)$,
where $\flat$ is given by
\begin{equation}
\nonumber
E_i^{\flat}=K_i^{-2}F_i, \; \; F_i^{\flat}=E_iK_i^2, 
\; \; K_i^{\flat}=K_i. 
\end{equation}

The $\sun$-module $*$-algebra structure on $\hc=(\h,\sharp)$
is given by
\begin{equation}
\nonumber
z_i^{\sharp}=\z_i, \; \; C^{\sharp}=C.
\end{equation}

Consider the real forms $\usu=(\usln,\natural)$ of $\usln$
parameterized by a sequence ${\boldsymbol\iota}=(\iota_0,
\iota_1,\dotsc,\iota_n)$, where $\iota_i=\pm 1$. These 
are given by
\begin{equation}
\nonumber
E_i^{\natural}=\iota_{i-1}\iota_iE_i^{\flat}, \; \;
F_i^{\natural}=\iota_{i-1}\iota_iF_i^{\flat}, \; \;
K_i^{\natural}=K_i.
\end{equation}
They are quantizations of different Lie bilagebra structures
on ${\frak g}_0={\frak su}(m,n+1-m)$, where $m$ is the number
of instances when $\iota_i=1$ and $n+1-m$ is the number of
instances when $\iota_i=-1$.

Then, the $\usu$-module $*$-algebra structure on $\ho=(\h,*)$
is given by
\begin{equation}
\label{*}
z_i^*=\iota_i\z_i, \; \; C^*=C.
\end{equation}
\end{exmp}

\begin{rem}
We will be particularly interested in the case when
${\boldsymbol \iota}=(-1,1)$. We denote $U_q{\frak su}
\langle -1,1\rangle$ by $\su$, $\ho$ by $\cal H$,
$E_1$ by $E$, $F_1$ by $F$, and $K_1$ by $K$.
\end{rem}

\section{Quantum Generalized Flag Manifolds}

In this section we establish the connection between the
quantum Heisenberg algebras and some quantum $G$-spaces.
From now on we use a shorter notation $\cal H$ for $\h$.

Recall a construction described in \cite{S}. Given a
finite-dimensional simple $\g$-nodule $V=L(\Lambda)$ with 
the highest weight $\Lambda$, we define a multiplication on
\begin{equation}
\nonumber
\CO^+=\bigoplus_{k=0}^{\infty}L(k\Lambda)
\end{equation}
as follows. Given $a\in L(k\Lambda)$ and $b\in L(m\Lambda)$,
we take the projection of $a\otimes b$ on $L((k+m)\Lambda)
\subset L(k\Lambda)\otimes L(m\Lambda)$ as the product of
$a$ and $b$ (note that it is correctly defined, as the multiplicity
of $L((k+m)\Lambda)$ in $L(k\Lambda)\otimes L(m\Lambda)$ is
equal to $1$). It is easy to see that it defines, in fact, a
$\g$-module algebra structure on $\CO^+$.

Define similarly a $\g$-module algebra by applying
the same construction to
\begin{equation}
\nonumber
\CO^-=\bigoplus_{k=0}^{\infty}L(k\Lambda)^*.
\end{equation}

\begin{PBWthm}
\label{PBWt}
The multiplication maps
\begin{gather}
\nonumber
\F^+\otimes\F^-\ri\F, \\
\nonumber
{\cal H}^+\otimes\h^0\otimes\h^-\ri\h
\end{gather}
are isomorphisms of $\g$-module algebras.
Here $\F^+={\cal H}^+$ is the subalgebra generated by
$V\subset V\oplus V^*$, $\F^-={\cal H}^-$ the subalgebra
generated by $V^*\subset V\oplus V^*$, and ${\cal H}^0$
the subalgebra generated by $C$.
\end{PBWthm}

Note that $\CO^+$ (resp. $\CO^-$) is the quantum analogue 
of the algebra of holomorphic (resp. anti-holomorphic) 
polynomial functions on the $G$-orbit of $\C v_{\Lambda}$ 
where $v_{\Lambda}\subset L(\Lambda)$ is a highest weight 
vector. At the same time, ${\cal H}^+=\F^+$ (resp. ${\cal H}^-=\F^-$) 
is a quantum analogue of the algebra of holomorphic (resp.
anti-holomorphic) polynomial functions on $V=L(\Lambda)$.

In the classical case (when $q=1$) $\C[{\cal O}_V]_1$ is 
a quotient of $\C[V]_1$ over the ideal generated by the
Pl\"{u}cker relations. A similar situation takes place
in the quantum case. Namely, it was shown in \cite{S} that 
$\CO^+$ is a quotient of ${\cal H}^+$ by an ideal $J_+$ generated
by the subspace 
\begin{equation}
\nonumber
E_{\Lambda}^+=\left(q^Z-q^{4(\Lambda+\rho,\Lambda)}\right)
\left(L(\Lambda)\otimes L(\Lambda)\right)
\end{equation} 
of quadratic relations called (holomorphic) {\em quantum 
Pl\"{u}cker relations\/}. Here $q^Z$ is the canonical central
element of $\g$ defined in \cite{D2}, and $\rho$ is the half 
of the sum of all positive roots of $\frak g$. Similarly, we 
can get that $\CO^-$ is a quotient of ${\cal H}^-$ over an ideal 
$J^-$ generated by the subspace
\begin{equation}
\nonumber
E_{\Lambda}^-=\left(q^Z-q^{4(\Lambda+\rho,\Lambda)}\right)
\left(L(\Lambda)^*\otimes L(\Lambda)^*\right)
\end{equation}
of what may be called anti-holomorphic quantum Pl\"{u}cker relations. 

Define the $\g$-module algebra $\CO$ as the quotient
of ${\cal H}$ over the ideal $J$ generated by all the quantum
Pl\"{u}cker relations in both $J^+$ and $J^-$. Instead 
of $C$ introduce a new generator
\begin{equation}
\nonumber
c=\frac{1}{q^{-1}-q}C.
\end{equation}

Now, if we take the quasi-classical limit $q\ra 1$, keeping $c$,
not $C$ constant, we will get a commutative Poisson algebra
which is the algebra of homogeneous polynomials on a family of
projective Poisson $G$-spaces with a homogeneous parameter $c$. 
Homogeneous, because $z_i$'s, $\z_i$'s and $c$ are defined
up to a group of automorphisms
\begin{equation}
\label{auto}
\kappa_{\alpha}:z_i\mra\alpha z_i, \; \; 
\kappa_{\alpha}:\z_i\mra\bar{\alpha}\z_i, 
\; \; \kappa_{\alpha}:c\mra |\alpha|^2c,
\end{equation}
for any $\alpha\in\C$. Note that the same formulas
define a group of automorphisms of the $\g$-module
algebra $\CO$ in the quantum case.

The above-mentioned projective $G$-spaces are the projectivizations
of the $G$-orbits of the space $\C v_{\Lambda}$ of highest weight 
vectors. They are called {\em generalized flag manifolds\/}. They 
are of the form $G/P$, where $P$ is the parabolic subgroup of $G$ 
whose Lie algebra is generated by ${\frak b}_+={\frak n}_+\oplus 
{\frak h}$ and the root vectors $E_{\lambda}^-$ such that $(\lambda,
\Lambda)=0$. In particular, if $\Lambda=\rho$, we get flag 
manifolds themselves. If $\Lambda$ is a fundamental weight, 
we get Grassmanians. This observation justifies the following
definition.

\begin{defn}
\label{flag}
The $\g$-module algebra $\CO$ is called the algebra of homogeneous 
polynomial functions on a family of quantum generalized flag manifolds.
\end{defn}

Let $G_c$ be the compact real form of $G$ whose quantization
yields $\gc$, $G_0$ the non-compact real form whose quantization
yields $\go$. It is easy to see that 
\begin{equation}
\nonumber
J^{\sharp}=J, \; \; J^*=J.
\end{equation}

Therefore, $\CO$ has canonical $\gc$- and $\go$-module
$*$-algebra structures. The $\gc$-module $*$-algebra
$(\CO,\sharp)$ can be thought of as the quantum algebra
of homogeneous polynomials on a family of generalized 
flag manifolds of the form $G/P$ considered as Poisson 
$G_c$-spaces. 

The $\go$-module $*$-algebra $(\CO,*)$ can be thought of 
as the quantum algebra of homogeneous polynomials on a 
family of corresponding symmetric Poisson $G_0$-spaces 
of non-compact type. As $G_0$-spaces they are isomorphic
to the $G_0$-orbit of the image of $P\subset G$ in $G/P$ 
with respect to the quotient map, where $P$ is the same
as above.

\begin{exmp}
\label{CO}
Return to Examples \ref{exmp-1} and \ref{exmp-2}. In this
case there are no Pl\"{u}cker relations, so that the
relations in $\CO$ look as follows:
\begin{eqnarray}
\nonumber
z_iz_j=qz_jz_i, \; (i<j), & \z_i\z_j=q^{-1}\z_j\z_i, \; (i<j), \\
\nonumber
z_i\z_j=q\z_jz_i, \; (i\not =j), & 
z_i\z_i-\z_iz_i=\left(q^{-2}-1\right)
\left(\sum_{k>i}z_k\z_k+qc\right), \\
\nonumber
z_i c=q^2c z_i, & \z_i c=q^{-2}c\z_i. 
\end{eqnarray}

The $\usu$-module $*$-algebra $(\CO,\sharp)$ is the
quantum algebra of homogeneous polynomials on a
family of quantum $\C P^n$, while the $\usu$-module
$*$-algebra $(\CO,*)$ is the quantum algebra of
homogeneous polynomials on a family of quantum
hyperboloids which posess a complex manifold
structure (inherited from $G/P$). 
\end{exmp}

\begin{rem}
\label{Podles}
When $G_c=SU(2)$, the family of quantum $\C P^1$
is nothing but the family of quantum Podle\'{s}
2-spheres introduced in \cite{P}.
\end{rem}

It is always nice to have a large commutative subalgebra. 
Let $U_q{\frak h}$ be the Hopf subalgebra of $\usln$
generated by $K_i, K_i^{-1}, \; i=1,2,\dotsc,n$. 
Consider the subalgebra $\Ch\subset\CO$ of the
$U_q{\frak h}$-invariant elements. Denote by
$\inv\subset\CO$ the subalgebra of $\g$-invariant
elements.

\begin{prop}
\label{Commute}
{\rom (1)} The algebra $\Ch$ is commutative and generated by
\begin{equation}
\nonumber
x_i=\sum_{k\geq i}z_k\z_k+qc, \; i=0,\dotsc,n+1. 
\end{equation}
Moreover, the following relations hold:
\begin{eqnarray}
\label{x1}
z_ix_j=q^2x_jz_i, \; (i<j), & \z_ix_j=q^{-2}x_j\z_i, \; (i<j), \\
\label{x2}
z_ix_j=x_jz_i, \; (i\geq j), & \z_ix_j=x_j\z_i, \; (i\geq j).
\end{eqnarray}

{\rom (2)} The algebra $\inv$ is generated by
\begin{equation}
\nonumber
c=q^{-1}x_{n+1} \; \; \text{and} \; \; d=qx_0.
\end{equation}
Moreover, $d$ belongs to the center of $\CO$.
\end{prop}

The formulas \eqref{x1}-\eqref{x2} allow to extend $\CO$
by adding functions of $\x=(x_0,x_1,\dotsc,x_{n+1})$ 
so that the following relations hold:
\begin{eqnarray}
\nonumber
z_if(\x) & = & f(x_0,\dotsc,x_i,q^2x_{i+1},
\dotsc,q^2x_{n+1})z_i, \\
\nonumber
\z_if(\x) & = & f(x_0,\dotsc,x_i,q^{-2}x_{i+1},
\dotsc,q^{-2}x_{n+1})\z_i.
\end{eqnarray}

We denote the extended algebra by $\Fun$. 

\begin{prop}
\nonumber
The $\gc$- and $\go$-module $*$-algebra structures can be
uniquely extended from $\CO$ on $\Fun$. The action of $\g$
is given by
\begin{eqnarray}
\nonumber
E_i:f(\x) & \mra & \frac{f(\x)-T_if(\x)}{x_i-q^2x_i}z_{i-1}\z_i,  \\
\nonumber
F_i:f(\x) & \mra & \frac{T_i^{-1}f(\x)-f(\x)}{q^{-2}x_i-x_i}z_i\z_{i-1},  \\
\nonumber
K_i:f(\x) & \mra & f(\x),
\end{eqnarray}
where
\begin{equation}
\nonumber
T_i:f(\x)\mra f(x_0,\dotsc,x_{i-1},q^2x_i,x_{i+1},\dotsc,x_{n+1}),
\end{equation}
with involutions given by 
\begin{equation}
\nonumber
f(\x)^{\natural}=f(\x), \; \; \;
f(\x)^*=f(\x).
\end{equation}
\end{prop}

\begin{exmp}
\label{SU11}
When $G_0=SU(1,1)$, one can modify our construction, 
so that we get more quantum $SU(1,1)$-spaces. 

Recall that $z_0,z_1,\z_0,\z_1$ are the homogeneous 
coordinates and $c$ a homogeneous parameter on a family 
of quantum $SU(1,1)$-spaces, and that they are defined
up to the automorphisms $\kappa_{\alpha}$ given by \eqref{auto}.

Consider a subalgebra $\FX\subset\Fun$ of $\kappa_{\alpha}$-invariant
elements in $\Fun$. Denote by $\widetilde{Func}(X)_q$ the subalgebra
in $\Fun$ of the elements which are invariant with respect to 
$\kappa_{\alpha}$ with $|\alpha|=1$.
\end{exmp}

\begin{prop}
\label{F-sl2}
The $\usl$-module algebra $\widetilde{Func}(X)_q$ is generated by
\begin{equation}
\nonumber
x=z_1\z_1+qc, \; \; y=z_0\z_1, \; \; \y=z_1\z_0
\end{equation}
with the relations
\begin{eqnarray}
\label{cd1}
yf(x)=f(q^2x)y, & \y f(x)=f(q^{-2}x)\y, \\
\label{cd2}
\y y=-(q^{-1}x-c)(q^{-1}x-d), & y\y=-(qx-c)(qx-d),
\end{eqnarray}
while $c$ and $d$ belong to the center of $\widetilde{Func}(X)_q$.
\end{prop}

It turns out that, besides the involution $*$ given in
\eqref{*}, there exists yet another one which makes
$\FX$ into a $\su$-module $*$-algebra. To keep the
notation shorter, we will use the somewhat larger
algebra $\widetilde{Func}(X)_q$.

\begin{prop}
\label{*-structures}
There are two non-equivalent $\su$-module $*$-algebra 
structures on $\widetilde{Func}(X)_q$, one of them given by
\begin{equation}
\label{su11*}
y^*=-\y, \; \; x^*=x, \; \; c^*=c, \; \; d^*=d, 
\end{equation}
and the other one given by
\begin{equation}
\label{star}
y^{\star}=-y, \; \; x^{\star}=x, \; \; c^{\star}=d.
\end{equation}
\end{prop}

In both cases we can define the $\su$-module $*$-algebra
$\cd$ as the quotient of $\widetilde{Func}(X)_q$ over the ideal generated
by $c-c_0$ and $d-d_0$, where $c_0,d_0\in{\Bbb R}$ in
the first case and $c_0\in\C, c_0=\bar{d_0}$ in the
secons case. 

It is clear that if $c_0,d_0\in{\Bbb R}$ and $c_0\not =d_0$,
$\cd$ is a quantum algebra of functions on the two-sheet
hyperboloid $|y|^2=(x-c_0)(x-d_0)$. If $c_0=d_0$, $\cd$ is a
quantum algebra of functions on the cone given by the same
equation. Finally, if $c_0=\bar{d_0}, c_0\not =d_0$, $\cd$
is a quantum algebra of functions on the corresponding 
one-sheet hyperboloid.

\section{Quantum Moment Map}

Recall the definition of the classical moment map, generalized
for the case when $G$ is a Poisson Lie group with a non-trivial
Poisson structure. Consider the corresponding Lie bialgebra
$\frak g$ and the dual Poisson Lie group $G^*$ which is defined 
as the connected and simply connected Poisson Lie group with the 
Lie bialgebra ${\frak g}^*$. For any $\xi\in {\frak g}$, let
$\alpha_{\xi}$ be the left invariant differential $1$-form on
$G^*$ with $\alpha_{\xi}(e)=\xi$. The Poisson bivector field
$\pi_{G^*}$ on $G$ defines a map $\check{\pi}_{G^*}:\Omega^1\left(
G^*\right)\rightarrow Vect\left(G^*\right)$.

A vector field $\rho_{\xi}=\check{\pi}_{G^*}(\alpha_{\xi})$
is called the left {\em dressing vector field\/} on $G^*$ corresponding
to $\xi\in {\frak g}$. The left dressing vector fields define a local
action of $G$ on $G^*$ which is called the left {\em dressing action\/}.
In some cases, for example, when $G$ is compact, it can be extended to
a global action. But in general, it need not be the case, as the example
of $G=SU(1,1)$ already shows.

Suppose now that $M$ is a left Poisson $G$-manifold,
that is, $M$ is a Poisson manifold with the action of $G$ on $M$
such that the corresponding map $G\times M\rightarrow M$ is a
Poisson map. Let $\sigma_{\xi}$ be the vector field corresponding
to the infinitesimal action of $\xi\in {\frak g}$. Keeping in mind 
that the local dressing action in our examples will not always be 
integrable to a global action, we modify slightly the usual 
definition of the generalized moment map (cf. \cite{Lu}) in 
order to apply it to our examples.

\begin{defn}
\label{cl-mm}
Let $M'$ be a union of symplectic leaves in $M$ such that $M'$ is
a dense subset in $M$. A map $J:M'\rightarrow G^*$ is called {\em 
a moment map\/} for $M$ if
\begin{equation}
\label{dr}
\sigma_{\xi}=\check{\pi}_M\left(J^*(\alpha_{\xi})\right).
\end{equation}
\end{defn}

We see that Definition \ref{cl-mm} means that
$J:M'\rightarrow G^*$ intertwines locally the $G$-action on $M$ 
with the dressing action of $G$ on $G^*$. When $G$ is a Poisson
ie group with the trivial Poisson structure, the dual Poisson
Lie group $G^*$ is isomorphic to ${\frak g}^*$ as a Poisson
manifold and is Abelian as a group. The corresponding dressing
action always extends to a global one which is nothing but the
usual coadjoint action of $G$ on ${\frak g}^*$. Thus, in this
case for any Hamiltonian $G$-space $M$, there exists a moment
map onto a coadjoin orbit in ${\frak g}^*$.

On the quantum level, it would have been reasonable to look for
a quantum moment map in the form $Func\left(G^*\right)_q\rightarrow 
Func(M)_q$. However, the Drinfeld's dualty tells us that the quantum
enveloping algebra $U_q{\frak g}$ can be thought of as a quantum
algebra of functions on $G^*$. Indeed, we will obtain a quantum
moment map in the form $U_q{\frak g}\mapsto Func(M)_q$.

As is well known, the quasi-classical analogue of the quantum
adjoin action of $U_q{\frak g}$ on itself given by
\begin{equation}
\label{q-adj}
ad_q(a):b\mapsto \sum_ka_k^{(1)}bS\left(a_k^{(2)}\right),
\; \; \mbox{\rm whenever} \; \; \Delta(a)=\sum_ka_k^{(1)}\otimes
b_k^{(2)} 
\end{equation}
is nothing but the left dressing action of $U{\frak g}$ on
$Func(G^*)$. Also, it is well known that for any Hopf algebra
$A$, the quantum adjoin action of $A$ on itself equips $A$ with
an $A$-module algebra structure, or an $A$-module $*$-algebra
structure if $A$ is a Hopf $*$-algebra.

This inspires the following definition, just slightly different
from the one given in \cite{Lu} (we do not have to worry about
the equality $M'=M$).

\begin{defn}
\label{q-mm}
\rom{(1)} Given a $U_q{\frak g}$-module algebra $\cal F$, a homomorphism
$J:U_q{\frak g}\rightarrow {\cal F}$ is called a {\em quantum
moment map\/} if $J$ is a morphism of $U_q{\frak g}$-module
algebras, with $U_q{\frak g}$ acting on itself by means of
the quantum adjoin action \eqref{q-adj}. 

\rom{(2)} Given a $U_q{\frak g}_0$-module $*$-algeba ${\cal F}_0$,
a $*$-homomorphism $J_0:U_q{\frak g}_0\rightarrow {\cal F}_0$ is called
a {\em quantum moment map\/} if $J_0$ is a morphism of $U_q{\frak g}_0$-module
$*$-algebras, with $U_q{\frak g}_0$ acting on itself by means of the quantum
adjoint action \eqref{q-adj}.
\end{defn}

We see that the quantum Heisenberg algebra ${\cal H}$ contains the subalgebras
${\cal H}^+$ and ${\cal H}^-$ generated
by $V$ and $V^*$ respectively. Of course, both are $U_q{\frak g}$-module
subalgebras of $\cal H$. Consider the subalgebra ${\cal H}^-_0$
generated by ${\cal H}^-$ and $C$. It has a one-dimensional
representation $\chi$ in ${\Bbb C}_{\chi}$ given by
\[ \chi(V)=0, \; \; \chi(C)=1. \]

Consider the corresponding induced $\cal H$-module
\[ W=Ind_{{\cal H}^-_0}^{\cal H}\left({\Bbb C}_{\chi}\right). \]

It is spanned by monomials of the form
\begin{equation}
\label{monomial}
a_{i_1}^{m_1}a_{i_2}^{m_2}...a_{i_k}^{m_k}{\boldsymbol 1}_{\chi},
\end{equation}
where $a_j\in V\subset {\cal H}^+$ and ${\boldsymbol 
1}_{\chi}$ is a generator of ${\Bbb C}_{\chi}$. Thus, we see that
$W$ is isomorphic to ${\cal H}^+$ as a vector space,
with a ${\Bbb Z}^{\mbox{\footnotesize dim}V}$-grading defined by 
\eqref{monomial}. This equips $W$ with a $U_q{\frak g}$-module 
structure so that $W$ is isomorphic to ${\cal H}^+$ as a $U_q{\frak g}$-module.

\begin{prop}
\label{repr}
\rom{(1)} The $\cal H$-module $W$ is simple and faithful.

\rom{(2)} The subalgebra ${\cal H}^{\rm inv}$ of
the $U_q{\frak h}$-invariant elements in $\cal H$ is commutative.
Moreover, any homogeneous monomial of the form \eqref{monomial}
in $W$ is an eigen-vector for the action of ${\cal H}^{\rm inv}$.
\end{prop}

The statement (1) of Proposition \ref{repr} shows that $\cal H$ is
isomorphic to its image in $End\,W$. Also, there exists a basis
in $W$ (spanned by the monomials of the form
\begin{equation}
\label{monimial}
f=v_0^{m_0}v_1^{m_1}...v_n^{m_n}{\boldsymbol 1}_{\chi}\in W
\end{equation}
which diagonalizes the action of ${\cal H}^{\rm inv}$.
This allows us to extend the algebra $\cal H$ by the functions on the
spectrum of ${\cal H}^{\rm inv}$ in $W$. Denote the
corresponding algebra by $\tilde{\cal H}$. One can show
that the $U_q{\frak g}$-module algera structure can be extended
from $\cal H$ to $\tilde{\cal H}$. 

Obviously, $\tilde{\cal H}$ is isomorphic to $End\,W$ 
as an algebra. On the other hand, $U_q{\frak g}$ acts in $W$.
This induces a homomorphism $J:U_q{\frak g}\rightarrow
\tilde{\cal H}$. It is clear that the image of 
$U_q{\frak g}$ lies in $Func(X)_q\subset\tilde{\cal H}$ -- 
the subalgebra of $\kappa_{\alpha}$-invariant elements in 
$\tilde{\cal H}$.

\begin{thm}
\label{J}
\rom{(1)} There exists a unique \rom{(}up to a $U_q{\frak g}$-module
algebra automorphism of $U_q{\frak g}$\rom{)} homomorphism 
\[ J:U_q{\frak g}\rightarrow Func(X)_q \]
such that the composition of $J$ with the action of
$Func(X)_q\subset\tilde{\cal H}$ in $W$ coincides 
with the action of $U_q{\frak g}$ in $W$.

\rom{(2)} $J$ is a morphism of $U_q{\frak g}$-module algebras,
with $U_q{\frak g}$ acting on itself via the quantum adjoin action.
In other words, $J$ is a quantum moment map for $Func(X)_q$.
\end{thm}

\begin{pf}
The first statement has been proved above. To show that (2) holds,
note that the image of $U_q{\frak g}$ in $\tilde{\cal H}$
must preserve the scalar degree $m=|{\boldsymbol m}|=m_0+m_1+...+m_n$
of a monomial of the form \eqref{monomial}. Therefore, $U_q{\frak g}$
maps into the subalgebra generated by the elements of the form $\varphi
\hat{\psi}$, where $\varphi$ (resp. $\hat{\psi}$) belongs to the 
subalgebra of $\kappa_{\alpha}$-invariant elements in the extension 
$\tilde{\cal H}^+$ (resp. $\tilde{\cal H}^-$) 
of ${\cal H}^+$ (resp. ${\cal H}^-$) by 
$\tilde{H}_{\cal g}(V)^{inv}$. In particular, $U_q{\frak h}$ can be 
shown to be mapped into the subalgebra of $\kappa_{\alpha}$-invariant
elements in $\tilde{\cal H}^{\rm inv}$.

Given $a\in U_q{\frak g}$ with $\Delta(a)=\sum_ka_k^{(1)}\otimes a_k^{(2)}$, 
we get that
\begin{eqnarray*} 
a\left(\varphi\hat{\psi}\right) & = & \sum_ka_k^{(1)}(\varphi)
\widehat{\left(S\left(a_k^{(2)}\right)\psi\right)}= \\
 & = & \sum_kJ\left(a_k^{(1)}\right)\varphi\widehat{\left(
J\left(S\left(a_k^{(2)}\right)\right)\psi\right)}.
\end{eqnarray*}

Recalling that $\hat{\psi}=\psi^{-1}f$, where $f\in\tilde{\cal 
H}^{\rm inv}$, one can prove (2) after some short computations.
\end{pf}

Note that the moment map $J$ for $Func(X)_q$ was obtained first as
a homomorphism from $U_q{\frak g}$ into the quantum Heisenberg algebra.
Thus, we see another manifestation of the Drinfeld's duality. In this
case, the same map has two quasi-classical analogues. One of them is
the homomorphism from the classical universal enveloping algebra to
the Heisenberg algebra which corresponds to a realization of $U{\frak g}$
by differential operators on a $G$-space. Another one is a moment map
for a family of generalized flag manifolds. Let us look at a few examples.

\begin{exmp}
\label{mm-1}
In the context of Example \ref{exmp-1} (that is, when ${\frak g}=
{\frak sl}(n+1)$ and $V$ being the first fundamental representation),
the map $J$ of Theorem \ref{J} is given by
\begin{eqnarray}
\label{J-sln-1}
J:E_i & \mapsto & \frac{\left(q^{-1}-q\right)^{\frac{1}{2}}}
{(x_{i-1}x_{i+1})^{\frac{1}{2}}}z_{i-1}\hat{z}_i, \\
\label{J-sln-2}
J:F_i & \mapsto & \frac{\left(q^{-1}-q\right)^{\frac{1}{2}}}
{x_i}z_i\hat{z}_{i-1}, \\
\label{J-sln-3}
J:K_i & \mapsto & \frac{x_i}{(x_{i-1}x_{i+1})^{\frac{1}{2}}}.
\end{eqnarray}

Moreover, given a Hopf $*$-algebra structure $\usu$ on $U_q{\frak sl}(n+1)$
and the involution \eqref{*} on $Func(X)_q$, we see that $J$ is in fact
a morphism of $\usu$-module $*$-algebras, thus defining a quantum moment
map for the $\usu$-module $*$-algebra $Func(X)_q$. 

Note that the subalgebra $Func(X)_q$ of the $\kappa_{\alpha}$-invariant
elements in $\tilde{H}$ is isomorphic to the quantum Weyl algebra
constructed in \cite{H}. Also, the quantum moment map $J$ is
equivalent to the quantum oscillator map from $U_q{\frak sl}(n+1,{\Bbb C})$
into the quantum Weyl algebra constructed there.
\end{exmp}

\begin{rem}
\label{mm-rem}
Of course, given any automorphism $I$ of the $U_q{\frak g}$-module algebra
structure on $U_q{\frak g}$ (defined by the quantum adjoint action), the
map $J\circ I:U_q{\frak g}_0\rightarrow Func(X)_q$ yields another quantum
moment map. In particular, in the above example the group of such 
automorphisms is generated modulo the center by the automorphisms given
by 
\[ I_i:E_j\mapsto (-1)^{\delta_{ij}}E_j, \; \; I_i:F_j\mapsto F_j,
\; \; I_i:K_j\mapsto (-1)^{\delta_{ij}}K_j. \]

It is easy to see that the corresponding moment maps $J_i=J\circ I_i$
are given by the same formulas \eqref{J-sln-1}-\eqref{J-sln-3} except
that we take another value of $(x_{i-1}x_{i+1})^{\frac{1}{2}}$.
\end{rem}

\begin{exmp}
\label{mm-2}
Consider the case ${\frak g}={\frak sl}(2)$. Recall that in this case
$Func(X)_q$ can be described in terms of the generators $x,y,\hat{y}$
as given in Proposition \ref{F-sl2}. Then we have
\begin{eqnarray}
\label{J-sl2-1}
J:E & \mapsto & \frac{\left(q^{-1}-q\right)^{\frac{1}{2}}}
{(cd)^{\frac{1}{2}}}y, \\
\label{J-sl2-2}
J:F & \mapsto & \frac{\left(q^{-1}-q\right)^{\frac{1}{2}}}{x}\hat{y}, \\
\label{J-sl2-3}
J:K & \mapsto & \frac{x}{(cd)^{\frac{1}{2}}}.
\end{eqnarray}

Moreover, given the Hopf $*$-algebra $U_q{\frak su}(1,1)$, $J$ is a
morpism of $U_q{\frak su}(1,1)$-module $*$-algebras for any of the
involutions $*$ and $\star$ on $Func(X)_q$ defined by \eqref{su11*}
and \eqref{star} respectively. Suppose that we fix $c=c_0$ and
$d=d_0$ so that $c_0d_0>0$. Then $J$ is a quantum moment map for 
any of the quantum hyperboloids (or a quantum cone) $X_{c_0,d_0}$ 
defined in the previous section. It is interesting to
note that the quantum quadratic Casimir element 
\begin{equation}
\label{C_q}
C_q=\frac{1}{2}(EF+FE)+\frac{q^{-1}+q}{2\left(q^{-1}-q\right)^2}
\left(K-2+K^{-1}\right) 
\end{equation}
is mapped to
\begin{equation}
\label{J-C}
J:C\mapsto \frac{1}{\left(q^{-1}-q\right)^2}\left(\frac{c_0}
{d_0}+\frac{d_0}{c_0}-q^{-1}-q\right). 
\end{equation}

The quasi-classical analogue $J_0$ of $J$ imbeds a dense subset of 
the hyperboloid (or cone) $|y|^2=(x-c_0)(x-d_0)$ (precisely, the one
defined by $x\not =0$) into the dual Poisson Lie group $SU(1,1)^*$.
The picture is as follows. One can show that $SU(1,1)^*$ is isomorphic
as a Lie group to the group of translations and dilations of a plane,
or to the group of the matrices of the form $\left(\begin{array}{cc}
t & z \\ 0 & t^{-1} \end{array}\right)$, where $t>0$ and $z\in {\Bbb C}$.
We can assume without loss of generality that $c_0d_0=1$. Then, the 
piece of the manifold $|y|^2=(x-c_0)(x-d_0)$ with $x>0$
maps into $|z|^2=(t-c_0)(t-d_0)$, while the piece with $x<0$ maps
into $|z|^2=(t+c_0)(t+d_0)$, which is nothing but the reflection
$t\mapsto -t$ of $|z|^2=(t-c_0)(t-d_0)$ with $t<0$. 

Of course, these imbeddings preserve the symplectic leaves. Indeed,
for the one-sheet hyperboloids, both pieces $x>0$ and $x<0$ are
two-dimensional symplectic leaves, while any point of the circle 
$x=0$ is a zero-dimensional symplectic leaf. For a twho-sheet
hyperboloid with $0<c_0<d_0$, the whole sheet $x\geq d_0$ is a
symplectic leaf, and the two-dimensional pieces $x<0$ and $0<x
\leq c_0$ of the other sheet are symplectic leaves as well, while
any point of the circle $x=0$ is a zero-dimensional symplectic leaf.
Finally, for the cone with $c_0=d_0=1$, the pieces with $x<0$, $0<x<1$,
and $x>1$ are symplectic leaves, as are the points on the circle
$x=0$ and the vertex of the cone -- the unit element of the group.
\end{exmp}

\section{Quantum Polarizations}

In the previous section we constructed a quantum moment map $J:U_q{\frak g}
\rightarrow Func(X)_q$. If we find now a way to construct an irreducible
$*$-representation $\pi$ of $Func(X)_q$, the composition $\pi\circ J\circ 
I$ (for some $U_q{\frak g}_0$-module $*$-algebra automorphism $I$ of 
$U_q{\frak g}_0$) will give us a $*$-representation of $U_q{\frak g}_0$. 
It will be irreducible, because the image of $J$ coincides with the 
subalgebra in $Func(X)_q$ of the elements invariant with respect to 
the automorphisms $\kappa_{\alpha}$.

Recall that the classical orbit method constructs an irreducible 
representation of an algebra of functions on a Hamiltonian manifold
in sections of a certain linear bundle with connection (whose curvature 
is equal to the symplectic form) which are constant along a given
polarization. Our construction in the quantum case draws the ideas
from that classical picture.

We have constructed $Func(X)_q$ as the subalgebra in $\tilde{H}$
which consists of the $\kappa_{\alpha}$-invariant elements. Recall that
$\kappa_{\alpha}$ is a family of automorphisms parameterized by a non-zero
complex number $\alpha$ (see \eqref{auto}). If we think of $Func(X)_q$ as 
an algebra of functions on a quantum space $X$, then $\tilde{H}$ can be 
thought of as an algebra of functions on the total space of a linear 
bundle over $X$.

Consider the subalgebra $Hol(X)_q^+$ (resp. $Hol(X)_q^-$) of 
in $Func(X)_q$ generated by $v_kv_m^{-1}$ (resp. $\hat{v}_k
\hat{v}_m^{-1}$), where $v_k$ are vectors of a $U_q{\frak h}$-invariant
basis in $V$. We will see that in the examples it is going to play
the role of the algebra of holomorphic (resp. anti-holomorphic)
functions in the case of a complex polarization. The following
proposition reflects the fact that we deal with
a quantum analog of a $G$-invariant polarization.

\begin{prop}
\label{G-inv}
$Hol(X)_q^{\pm}$ is a $U_q{\frak g}$-module subalgebra in
$Func(X)_q$.
\end{prop}

\begin{exmp}
\label{hol-sln}
Suppose that, as in Example \ref{exmp-1}, ${\frak g}={\frak sl}(n+1,
{\Bbb C})$ (equipped with the standard Lie bialgebra structure), and
$V$ is the highest weight $U_q{\frak g}$-module with the highest weight
$\omega_1$ (the first fundamental weight). We keep the same notation
as before.

Then $Hol(X)_q^+$ is generated by 
\[ \zeta_i=z_i^{-1}z_{i-1}, \; \; \mbox{\rm ($i=1,2,...,n$)} \]
with the relations
\begin{eqnarray}
\label{zz-1}
\zeta_i\zeta_j & = & q^{\pm 1}\zeta_j\zeta_i, \; \mbox{\rm if} \; j=i\pm 1, \\
\label{zz-2}
\zeta_i\zeta_{i+1} & = & \zeta_{i+1}\zeta_i, \l \mbox{\rm otherwise}. 
\end{eqnarray}
Respectively, $Hol(X)_q^-$ is generated by
\[ \hat{\zeta}_i=\z_{i-1}\z_i^{-1}, \; \; \mbox{\rm ($i=1,2,...,n$)} \]
with the relations
\begin{eqnarray}
\label{zz-3}
\hat{\zeta}_i\hat{\zeta}_j & = & q^{\mp 1}\hat{\zeta}_j\hat{\zeta}_i,
 \; \mbox{\rm if} \; j=i\pm 1, \\
\label{zz-4}
\hat{\zeta}_i\hat{\zeta}_{i+1} & = & \hat{\zeta}_{i+1}\hat{\zeta}_i,
\; \mbox{\rm otherwise}. 
\end{eqnarray}

It is easy to write down explicit formulas for the action of
$U_q{\frak sl}(n+1,{\Bbb C})$ in $Hol(X)_q^{\pm}$, but we will
do it only in the special case of $n=1$ (see below).
\end{exmp}

\begin{exmp}
\label{hol-sl2}
When $n=1$, $Hol(X)_q^+$ is generated by a single element 
\begin{equation}
\label{zeta-defn}
\zeta=z_1^{-1}z_0=(qx-c_0)^{-1}y, 
\end{equation}
while $Hol(X)_q^-$ is generated by
\[ \hat{\zeta}=\z_0\z_1^{-1}=\hat{y}(qx-c_0)^{-1}. \]

The action of $U_q{\frak sl}(2,{\Bbb C})$ in $Hol(X)_q^+$
is given by
\begin{eqnarray*}
E:f(\zeta) & \mapsto & -q\zeta^2\frac{f(\zeta)-f\left(\zeta q^2\right)}
{\zeta-\zeta q^2}, \\
F:f(\zeta) & \mapsto & \frac{f\left(\zeta q^{-2}\right)-f(\zeta)}
{\zeta q^{-2}-\zeta}, \\
K:f(\zeta) & \mapsto & f\left(\zeta q^{-2}\right).
\end{eqnarray*}

The action of $U_q{\frak sl}(2,{\Bbb C})$ in $Hol(X)_q^-$ is given
by similar formulas. One can check that the center of $U_q{\frak sl}
(2,{\Bbb C})$ acts trivially in $Hol(X)_q^{\pm}$.
\end{exmp}

\begin{prop}
\label{zzx}
The $U_q{\frak sl}(n+1,{\Bbb C})$-module algebra $Func(X)_q$ is generated 
by $\zeta_i$, $\hat{\zeta}_i$ \rom{(}$i=1,2,...,n$\rom{)}, and the functions 
$f(x_1,x_2,...,x_n)$ with the relations \eqref{zz-1}-\eqref{zz-4} and
\begin{eqnarray}
\label{zx-1}
\zeta_if(x_1,...,x_n) & = & 
f\left(x_1,...,x_{i-1},q^2x_i,x_{i+1},...,x_n\right)\zeta_i, \\
\label{zx-2}
\hat{\zeta}_if(x_1,...,x_n) & = & 
f\left(x_1,...,x_{i-1},q^{-2}x_i,x_{i+1},...,x_n\right)\hat{\zeta}_i, \\
\label{zx-3}
\zeta_i\hat{\zeta}_i & = & \frac{x_{i-1}-x_i}{x_i-q^{-2}x_{i+1}}, \\
\label{zx-4}
\hat{\zeta}_i\zeta_i & = & \frac{q^2x_{i-1}-x_i}{x_i-x_{i+1}}.
\end{eqnarray}

Recall that given a real form $\usu$ of $U_q{\frak sl}(n+1,{\Bbb C})$, 
the involution \eqref{*} equips $Func(X)_q$ with a $\usu$-module $*$-algebra
structure. The involution \eqref{*} is given in the above generators by
\[ \zeta_i^*=\iota_i\hat{\zeta}_i, \; \; x_i^*=x_i. \]
\end{prop}

In the next section we will consider first the simplest case of 
$U_q{\frak su}(1,1)$ to illustrate the basic ideas. We will use 
them later to construct some irreducible $*$-representations of 
$\usu$ which correspond to the dressing orbits of the minimal 
dimension. Let us describe, therefore, the relations in that
special case more explicitly.

\begin{cor}
\label{zzx2}
\rom{(1)} For $c_0d_0=1$, $c_0,d_0\in {\Bbb R}$, the $U_q{\frak su}
(1,1)$-module $*$-algebra $Func(X_{c_0,d_0})_q$ is generated by $\zeta$, 
$\zeta^*$, and the functions $f(x)$ with the relations
\begin{eqnarray}
\label{*-z1}
\zeta f(x)=f\left(q^2x\right)\zeta, & 
\zeta^*f(x)=f\left(q^{-2}x\right)\zeta^*, \\
\label{*-z2}
\zeta\zeta^*=\frac{x-q^{-1}d_0}{x-q^{-1}c_0}, &
\zeta^*\zeta=\frac{x-qd_0}{x-qc_0}.
\end{eqnarray}

In particular, the following relation holds:
\[ \Phi\left(\zeta\zeta^*\right)=q^2\Phi\left(\zeta^*\zeta\right), \]
where
\begin{equation}
\label{Phi}
\Phi(t)=\frac{1-\gamma t}{1-t}, \; \; \gamma=\frac{c_0}{d_0}.
\end{equation}

\rom{(2)} For $c_0d_0=1$, $c_0,d_0\not\in {\Bbb R}$, the $U_q{\frak su}
(1,1)$-module $*$-algebra $Func(X_{c_0,d_0})_q$ is generated by $\zeta$,
$\hat{\zeta}$, and the functions $f(x)$ with the relations
\begin{eqnarray}
\label{star-z1}
\zeta f(x)=f\left(q^2x\right)\zeta, & 
\hat{\zeta}f(x)=f\left(q^{-2}x\right)\hat{\zeta}, \\
\label{star-z2}
\zeta\hat{\zeta}=\frac{x-q^{-1}d_0}{x-q^{-1}c_0}, &
\hat{\zeta}\zeta=\frac{x-qd_0}{x-qc_0}, \\
\label{star-z3}
\zeta\zeta^{\star}=\zeta^{\star}\zeta=1, & 
\hat{\zeta}\hat{\zeta}^{\star}=\hat{\zeta}^{\star}\hat{\zeta}=1.
\end{eqnarray}

In particular, the following relation holds:
\[ \Phi\left(\zeta\hat{\zeta}\right)=q^2\Phi\left(\hat{\zeta}\zeta\right), \]
where $\Phi(t)$ is given by \eqref{Phi}.
\end{cor}

The algebras $Hol(X)_q^{\pm}$ play the role of the
algebras of functions which are constant along a polarization. In
particular, in the cases of $n>1$ and of $c_0,d_0\in {\Bbb R}$ ($n=1$)
the corresponding polarization is complex, so that they are quantum
analogs of the algebras of holomorphic and anti-holomorphic functions.
For the quantum one-sheet hyperboloids $c_0,d_0\not\in {\Bbb R}$
($n=1$), however, the polarization is real, and $Hol(X_{c_0,d_0})_q^{\pm}$
are quantum analogs of the algebras of functions which are constant
along the two families of straight lines on the corresponding 
one-sheet hyperboloids.

\section{Irreducible $*$-Representations Of $U_q{\frak su}(1,1)$}

In this section we will construct irreducible $*$-representations of
$U_q{\frak su}(1,1)$ associated with the quantum spaces $X_{c_0,d_0}$.
We assume that $c_0d_0=1$.

As follows from \eqref{*-z2} and \eqref{star-z2}, $Func(X_{c_0,d_0})_q$
is generated as algebra by $\zeta$, $\zeta^{-1}$, and the functions
$f(x)$ (resp. by $\zeta^*$, $\left(\zeta^*\right)^{-1}$, and the functions
$f(x)$, or by $\hat{\zeta}$, $\hat{\zeta}^{-1}$, and $f(x)$), since we
see that, for instance, $\zeta^*=\frac{x-qd_0}{x-qc_0}\zeta^{-1}$.
In particular, we see that $Func(X_{c_0,d_0})_q$ is generated by two
distinguished subalgebras -- $Hol(X_{c_0,d_0})_q^{\pm}$ and $Func(
X_{c_0,d_0})_q^{inv}$.

In the classical case the orbit method realizes representations in
sections of a linear bundle with connection whose curvature coincides 
with the symplectic form on the corresponding coadjoint orbit. This is
a general fact that, given a symplectic form, if there exist any linear 
bundles with connection and that form as its curvature, they are
parameterized by the local systems on the manifold. Below we give
a construction of representations of $U_q{\frak su}(1,1)$. It will
be shown in the section about the quasi-classical analogs that the
character of the commutative subalgebra $Func(X_{c_0,d_0})_q^{inv}$
-- the one generated by the functions $f(x)$ -- plays the role of
a local system on the corresponding symplectic leaf.

Consider a $*$-homomorphism $\nu:Func(X_{c_0,d_0})_q^{inv}\rightarrow 
{\Bbb C}$ of the form
\[ \nu:f(x)\mapsto f(\nu_0), \; \mbox{where} \; 
\nu_0\in {\Bbb R}\setminus\{0\}. \]

It defines a one-dimensional $Func(X_{c_0,d_0})_q^{inv}$-module
${\Bbb C}_{\nu}$. Consider the induced $Func(X_{c_0,d_0})_q$-module
\[ \Pi_{\nu}=Ind_{F^{inv}}^F{\Bbb C}_{\nu}, \]
where we use a short notation $F=Func(X_{c_0,d_0})_q$.

\begin{prop}
\label{pi-cd}
The $Func(X_{c_0,d_0})_q$-module $\Pi_{\nu}$ is isomorphic to
$Hol(X_{c_0,d_0})_q^{\pm}$ as a $Hol(X_{c_0,d_0})_q^{\pm}$-module
\rom{(}with respect to the left multiplications\rom{)}. Moreover,
the action of $Func(X_{c_0,d_0})_q$ in $\Pi_{\nu}$ is given by
\begin{eqnarray}
\label{action-z1}
\zeta:f(\zeta){\boldsymbol 1}_{\nu} & \mapsto & 
\zeta f(\zeta){\boldsymbol 1}_{\nu}, \\
\label{action-z2}
x:f(\zeta){\boldsymbol 1}_{\nu} & \mapsto &
\nu_0f\left(\zeta q^{-2}\right){\boldsymbol 1}_{\nu}, \\
\label{action-z3}
y:f(\zeta){\boldsymbol 1}_{\nu} & \mapsto &
\zeta^2\frac{q^{-1}\nu_0f\left(\zeta q^{-2}\right)-c_0f(\zeta)}
{\zeta}{\boldsymbol 1}_{\nu}, \\
\label{action-z4}
\hat{y}:f(\zeta){\boldsymbol 1}_{\nu} & \mapsto &
-\frac{q\nu_0f\left(\zeta q^{-2}\right)-d_0f(\zeta)}
{\zeta}{\boldsymbol 1}_{\nu},
\end{eqnarray}
when $\Pi_{\nu}$ is realized as the span of monomials of the form
$\zeta^k{\boldsymbol 1}_{\nu}$, where ${\boldsymbol 1}_{\nu}$ is
a generator of ${\Bbb C}_{\nu}$. Similar formulas hold in the case
if we realize $\Pi_{\nu}$ as the span of monomials of the form
$\left(\zeta^*\right)^k{\boldsymbol 1}_{\nu}$ or
$\hat{\zeta}^k{\boldsymbol 1}_{\nu}$.
\end{prop}

As we see from \eqref{action2}, the set of eigen-values of the 
action of $x$ in $\Pi_{\nu}$ is a part of the geometric progression 
\[ {\frak M}_{\nu_0}=\{\nu_0q^{2k}\}_{k\in {\Bbb Z}}. \]

\begin{defn}
\label{unit-mod}
Suppose that $F$ is a $*$-algebra.

(1) We call a $F$-module $\Pi$ {\em unitarizable\/} 
if there exists a positive definite Hermitian scalar product 
$( \; , \; )$ in $\Pi$ such that
\[ (av_1,v_2)=(v_1,a^*v_2) \]
for any $a\in F$ and $v_1,v_2\in \Pi$.

(2) Suppose that $\Pi$ is a unitarizable $F$-module. Consider the
Hibert space $H$ which is the completion of $\Pi$. We say that the
action of $F$ in $\Pi$ defines a $*$-{\em representation\/} $\pi$ of
$F$ in $H$ if the action of any element $a\in F$ in $\Pi$ can be
extended to a closed operator $\pi(a)$ in $H$.
\end{defn}

\begin{thm}
\label{sp-x}
\rom{(1)} Let $0<c_0\leq d_0$ \rom{(}the case of the quantum two-sheet
hyperboloids and the quantum cone\rom{)}. Suppose that neither $x_2=qc_0$ 
nor $x_0=q^{-1}d_0$ belongs to ${\frak M}_{\nu_0}$. Then there exists a 
scalar product $( \; ,\; )$ in $\Pi_{\nu}$ making in into a simple 
unitarizable $Func(X_{c_0,d_0})_q$-module if and only if no point of 
${\frak M}_{\nu_0}$ lies in the interval $(qc_0,qd_0)$. The corresponding 
scalar product in $\Pi_{\nu}$ can be given by
\begin{equation}
\label{scalar}
(f(\zeta){\boldsymbol 1}_{\nu},g(\zeta){\boldsymbol 1}_{\nu})=
\nu\left(g(\zeta)^*f(\zeta)\right), 
\end{equation}
where $\nu$ is extended to $Func(X_{c_0,d_0})_q$ by $\nu(\zeta)=
\nu\left(\zeta^*\right)=0$. Moreover, the action of $Func(X_{c_0,d_0})_q$
in $\Pi_{\nu}$ defines an irreducible $*$-representation of
$Func(X_{c_0,d_0})_q$. The spectrum of the action of $x$ in the
corresponding Hilbert space is equal to ${\frak M}_{\nu_0}\cup\{0\}$.

\rom{(2)} Let $c_0=\bar{d}_0\not\in {\Bbb R}$ \rom{(}the case of the
quantum one-sheet hyperboloids\rom{)}. Then there exists a scalar product
$( \; , \; )$ in $\Pi_{\nu}$ making it into a unitarizable $Func(X_{c_0,
d_0})_q$-module. It can be given by \eqref{scalar}. Moreover, the action
of $Func(X_{c_0,d_0})_q$ in $\Pi_{\nu}$ defines an irreducible 
$*$-representation of $Func(X_{c_0,d_0})_q$. The spectrum of the 
action of $x$ in the corresponding Hilbert space is equal to 
${\frak M}_{\nu_0}\cup\{0\}$.

\rom{(3)} Let $0<c_0\leq d_0$. Suppose that $qc_0\in {\frak M}_{\nu_0}$
\rom{(}resp. $q^{-1}d_0\in {\frak M}_{\nu_0}$\rom{)}. Then there exists
a scalar product
$( \; , \; )$ in $\Pi_{\nu}$ making it into a unitarizable $Func(X_{c_0,
d_0})_q$-module. It can be given by \eqref{scalar}. Moreover, the action
of $Func(X_{c_0,d_0})_q$ in $\Pi_{\nu}$ defines an irreducible 
$*$-representation of $Func(X_{c_0,d_0})_q$. The spectrum of the 
action of $x$ in the corresponding Hilbert space is equal to 
\[ {\frak M}_+=\{c_0q^{2k+1}\}_{k=0}^{\infty}\cup\{0\} \; \; 
\mbox{({\em resp.\/}} \; {\frak M}_-=\{d_0q^{-2k-1}\}_{k=0}^{\infty}
\mbox{)}. \]
\end{thm}

\begin{pf}
It is easy to see that any monomial
of the form $\zeta^k{\boldsymbol 1}_{\nu}\in\Pi_{\nu}$ is an eigen-vector
for the action of $x$ with the eigen-value $\nu_0q^{-2k}$. At the same
time \eqref{action-z1}-\eqref{action-z4} show that the set of eigen-values
of $x$, being a part of the geometric progression ${\frak M}_{\nu_0}$
would truncate only if either $qc_0$ or $q^{-1}d_0$ belong to 
${\frak M}_{\nu_0}$. This follows also from \eqref{zeta-defn}.
Therefore, we need to show only the unitarizability of $\Pi_{\nu}$.
We need to remind some definitions.

\begin{defn}
\label{integral}
Suppose that $A$ is a Hopf $*$=algebra, $F$ an $A$-module $*$-algebra.
A linear functional $f\mapsto \int fd\mu$ defined on a linear subset $F_0$
of $F$ is called an {\em invariant integral\/} on $F$ if the following
properties are satisfied:
\begin{eqnarray*}
\int afd\mu=\varepsilon(a)\int fd\mu, \; \;
\int f^*d\mu=\overline{\int fd\mu}, \\
f\mapsto \int f^*fd\mu \; \mbox{is a positive definite form on $F_0$,}
\end{eqnarray*}
for any $a\in A$ and $f\in F_0$, where $\varepsilon$ is the counit in $A$.
\end{defn}

The following lemma is well known.

\begin{lem}
\label{obvious}
Suppose that $A$ is a Hopf $*$-algebra, $F$ an $A$-module $*$-algebra
with an invariant integral $\int d\mu:F_0\rightarrow {\Bbb C}$. Consider
the space $L_2(F,d\mu)$ consisting of all $f\in F$ such that $\int f^*fd\mu
<\infty$. Then $L_2(F,d\mu)$ is a Hilbert space with the scalar product
given by
\begin{equation}
\label{int-trivial}
(f,g)=\int g^*fd\mu,
\end{equation}
and the action of $F$ in $L^2(F,d\mu)$ by left multiplication defines
a $*$-representation of $F$.
\end{lem}

\begin{prop}
\label{inv-int}
Under the assumptions of Theorem \ref{sp-x} \rom{(1)-(2)}, the linear 
functional
\begin{equation}
\label{int-1}
\int\zeta^kf(x)d\mu=\delta_{k,0}\left(q^{-1}-q\right)
\sum_{x\in {\frak M}_{\nu_0}}xf(x) 
\end{equation}
is an invariant integral on $Func(X_{c_0,d_0})_q$. Similarly, under
the assumptions of Theorem \ref{sp-x} \rom{(3)}, the linear functional
\begin{equation}
\label{int-2}
\int\zeta^kf(x)d\mu=\delta_{k,0}\left(q^{-1}-q\right)
\sum_{x\in {\frak M}_{\pm}}xf(x)
\end{equation}
is an invariant integral on $Func(X_{c_0,d_0})_q$.
\end{prop}

\begin{pf}
This follows from the fact that $J\left(q^{\rho}\right)=J(K)=x$, where
$\rho\in U_q{\frak g}$ is the half of the sum of all positive roots.
But it can also be checked by a straightforward computation.
\end{pf}

Now we can prove Theorem \ref{sp-x}. We realize $\Pi_{\nu}$
as a subspace in $Func(X_{c_0,d_0})_q$ by mapping ${\boldsymbol 1}_{\nu}$
into the function $\delta_{\nu_0}$ which takes the value $1$ at $\nu_0\in 
{\frak M}_{\nu_0}$ and $0$ at any other point of ${\frak M}_{\nu_0}$.
It is easy to see that this map intertwines the action of $Func(X_{c_0,d_0})_q$
in $\Pi_{\nu}$ and itself by left multiplications. By Lemma \ref{obvious},
\eqref{int-trivial} defines a $*$-representation of $Func(X_{c_0,d_0})_q$
in $\Pi_{\nu}$ when $\int d\mu$ is of the form \eqref{int-1}-\eqref{int-2}.
Obviously, this representation is irreducible, since the module $\Pi_{\nu}$
is simple. On the other hand, \eqref{*-z2} forbids any point from the
spectrum of $x$ to lie in the interval $(qc_0,qd_0)$.
\end{pf}

The full list of irreducible $*$-representations of $U_q{\frak su}(1,1)$
was described in a number of papers (cf. \cite{Jap,KVrep}). They are
parameterized by the spin $l\in {\Bbb C}$ and the parity $|\epsilon|\leq
\frac{1}{2}$. The corresponding representations $T_{l,\epsilon}$ are 
subjects to the symmetries
\[ T_{l,\epsilon}\simeq T_{-l-1,\epsilon}, \; \;
T_{l,\epsilon}\simeq T_{l+\frac{2\pi i}{\log q},\epsilon}. \]

They are divided into the following series:
\begin{enumerate}
\item principal continuous series: $l=-\frac{1}{2}+i\rho$,
\item complimentary series: $l\in (-\frac{1}{2},0)$, $|\epsilon|<|l|$,
\item holomorphic discrete series: $l-\epsilon\in {\Bbb Z}$,
\item anti-holomorphic discrete series: $l+\epsilon\in {\Bbb Z}$,
\item strange series: $\mbox{Im}\,l=\frac{\pi i}{\log q}$.
\end{enumerate}

The can be distinguished by the action of the quantum quadratic Casimir
element \eqref{C_q} as $\frac{\left(q^l-q^{-l}\right)\left(q^{l+1}-
q^{-l-1}\right)}{\left(q^{-1}-q\right)^2}$ and by the fact that any 
eigen-value of $K$ in $T_{l,\epsilon}$ is of the form $q^{k+\epsilon}$.

Note that the strange series does not survive in the classical
limit. We will see later how it manifests itself in the language of the
symplectic leaves in $SU(1,1)^*$.

\begin{thm}
\label{realization}
Given an irreducible $*$-representation $\pi$ of $FUnc(X_{c_0,d_0})_q$,
$\pi\circ J$ is an irreducible $*$-representation of $U_q{\frak su}(1,1)$.

\rom{(1)} The irreducible $*$-representations described in Theorem \ref{sp-x}
\rom{(1)} give rise to the complimentary series representations if $\nu_0>0$
and the strange series representations if $\nu_0<0$.

\rom{(2)} The irreducible $*$-representations described in Theorem \ref{sp-x}
\rom{(2)} give rise to the principal continuous series representations.

\rom{(3)} The irreducible $*$-representations described in Theorem \ref{sp-x}
\rom{(3)} give rise to the holomorphic discrete series representations if 
they correspond to ${\frak M}_+$ and the anti-holomorphic discrete series 
representations if they correspond to ${\frak M}_-$.
\end{thm}

\begin{pf}
The theorem follows immediately from \eqref{J-C} if we assume $c_0=
q^{l+\frac{1}{2}}$ and $d_0=q^{-l-\frac{1}{2}}$ and from \eqref{J-sl2-3},
since any eigen-value of $K$ will be an eigen-value of $x$.
\end{pf}

Note that the principal continuous series representations correspond to the
symplectic leaves that are the halves ($x>0$ and $x<0$) of the one-sheet 
hyperboloids $X_{c_0,d_0}$ as described in Example \ref{mm-2}. The 
holomorphic discrete series representations and the strange series 
representations correspond to the halves ($x<0$ and $0<x\leq c_0$) 
of a sheet of the corresponding two-sheet hyperboloids $|y|^2=(x-c_0)
(x-d_0)$, while the anti-holomorphic series representations correspond
to the other sheet of the two-sheet hyperboloid. 

What is especially interesting is that the complimentary series 
representations correspond to the case when a geometric progression
${\frak M}_{\nu_0}$ can jump over the narrow interval $(qc_0,qd_0)$.
It looks as if in the quantum case the invariant measure can be extended
from one sheet of a quantum two-sheet hyperboloid onto another one, thus
making it 'connected'. We will call such quantum hyperboloids {\em quantum
tunnel hyperboloids\/}. This effect disappears in the classical limit.
In particular, it is interesting to compare it with the fact that
the classical orbit method fails to realize the complimentary series
representations.

As we pass to the quasi-classical limit, we see that the choice of
a geometric progression reflects the value of the parity $\epsilon$
of the corresponding irreducible $*$-representation of $U_q{\frak su}(1,1)$.
Therefore, we can think of the choice of ${\Bbb C}_{\nu}$ as the choice of 
a local system on the corresponding symplectic leaf in $SU(1,1)^*$.
However, with the 'tunnel effect' in mind, we see that certain choices
of ${\Bbb C}_{\nu}$ may not have any local system as their quasi-classical
analogs.

On the other hand, the observed correspondence between the symplectic 
leaves in $SU(1,1)^*$ and the representations of $U_q{\frak su}(1,1)$
depends on $q$, as $\frac{c_0}{d_0}=q^{2l+1}$. As we keep the spin $l$
of the representation fixed and take the limit $q\rightarrow 1$, the
corresponding symplectic leaves face two options. Those which give
rise to the strange series representations will go to infinity (thus,
there are no strange series representations in the classical case).
The other ones will tend to the nilpotent cone $|y|^2=(x-1)^2$.
If we consider them as points in the corresponding orbifold, we can
look at the rate with which the corresponding curve in the orbifold
tends to the cone. It will be an orbit of the coadjoint action
in ${\frak su}(1,1)^*$. 

This gives us the usual correspondence between the representations and 
coadjoint orbits described by the classicaal orbit method. Except that 
the 'tunnel effect' will disappear, and so will the above geometric 
realization of the complimentary series representations.

\section{Degenerate Series Of Irreducible $*$-Representations}

Consider the quantum moment map $J:U_q{\frak g}\rightarrow Func(X)_q$.
If we find now a way to construct an irreducible $*$-representation 
$\pi$ of $Func(X)_q$, the composition $\pi\circ J\circ I$ (for some 
$U_q{\frak g}_0$-module $*$-algebra automorphism $I$ of $U_q{\frak g}_0$) 
will give us a $*$-representation of $U_q{\frak g}_0$. It will be irreducible, 
because the image of $J$ coincides with the subalgebra in $Func(X)_q$ 
of the elements invariant with respect to the automorphisms $\kappa_{\alpha}$.

Throughout the section we consider the case $U_q{\frak g}_0=\usu$.
We proceed in a similar way as we did when $U_q{\frak g}_0=U_q{\frak su}
(1,1)$. Namely, consider a one-dimensional $*$-representation $\chi:
Func(X)_q^{inv}\rightarrow {\Bbb C}_{\chi}$ of $Func(X)_q^{inv}$. Recall
that $Func(X)_q^{inv}$ is commutative and generated by $x_0=q^{-1}d,
x_1,...,x_n,x_{n+1}=qc$. Since $x_i^*=x_i$, we can think of $\chi$ as
a triple $(c_0,d_0,\hat{\chi})$, where $c_0=\chi(c)$, $d_0=\chi(d)$, 
and $\hat{\chi}=(\chi(x_1),...,\chi(x_n))$ the point in ${\Bbb R}^n$.

Since $x_1,...,x_n$ are invertible in $Func(X)_q$, $\chi(x_i)\not =0$
for any $i=1,...,n$. If $\chi(x_i)<0$, we can choose the $\usu$-module
$*$-algebra automorphism $I$ so that, after replacing $\chi$ by $\chi\circ I$,
it becomes positive (see Remark \ref{mm-rem}). Thus, without losing
generality, we may assume that $\chi(x_i)>0$ for any $i=1,...,n$.

Consider an isomorphism of vector spaces ${\Bbb R}^n\rightarrow {\frak h}^*$
such that 
\[ (\lambda_1,...,\lambda_n)\mapsto\lambda_1\alpha_1+...+\lambda_n\alpha_n, \]
where $\alpha_1,...,\alpha_n$ are the positive simple roots. Thus,
we can think of $\hat{\chi}$ as a point in $q^{{\frak h}^*}$, so that
$\hat{\chi}=q^{\alpha}$, where $\alpha\in {\frak h}^*$.

\begin{prop}
\label{Gen-1}
The induced $Func(X)_q$-module
\[ W=Ind_{Func(X)_q^{inv}}^{Func(X)_q}{\Bbb C}_{\chi}. \]
is isomorphic to $Hol(X)_q^{\pm}$ as a $Hol(X)_q^{\pm}$-module.
\end{prop}

The following proposition follows immediately from Proposition \ref{Gen-1} 
and \eqref{zx-1}-\eqref{zx-2}.

\begin{prop}
\label{Gen-2}
$W$ is spanned by the common eigen-vectors of $x_1,...,x_n$, the set
$\cal Q$ of the corresponding eigen-values forming a part of the lattice 
$q^{P+\alpha}\subset q^{{\frak h}^*}$, where $P\subset {\frak h}^*$ is 
the weight lattice.
\end{prop}

Let $v\in W$ be such an eigen-vactor, with the eigen-value $(\lambda_1,
...\lambda_n)$. By \eqref{zx-3}-\eqref{zx-4}, it follows that
\begin{eqnarray}
\label{zl-1}
\zeta_i\zeta_i^*:v & \mapsto & \iota_i\frac{\lambda_{i-1}-\lambda_i}
{\lambda_i-q^{-2}\lambda_{i+1}}v, \\
\label{zl-2}
\zeta_i^*\zeta_i:v & \mapsto & \iota_i\frac{q^2\lambda_{i-1}-\lambda_i}
{\lambda_i-\lambda_{i+1}}v.
\end{eqnarray}

\begin{prop}
\label{Gen-3}
The linear functional
\begin{equation}
\label{int-rho}
\nu_W(f)=\mbox{\rm tr}_W\left(f(J\circ I)\left(q^{\rho}\right)\right),
\end{equation}
where $\rho$ is half the sum of all positive roots, is an invariant
integral on $Func(X)_q$.
\end{prop}

\begin{pf}
This easily follows from the fact that $J\circ I$ is a quantum moment
map and, thus, intertwines the $\usu$-action on $Func(X)_q$ with its
quantum adjoint action on itself, and from the well-known properties
of the distinguished element $\rho$.
\end{pf}

\begin{thm}
\label{Gen-4}
Let $\pi$ be the representation of $Func(X)_q$ in $W$, and a $\usu$-module
$*$-algebra automorphism $I$ of $\usu$ is chosen so that $\chi(x_i)>0$ for
any $i=1,...,n$. Then $\pi\circ J\circ I$ is an irreducible $*$-representation
of $\usu$ if and only if 
\begin{eqnarray}
\label{ll-1}
\iota_i\frac{\lambda_{i-1}-\lambda_i}
{\lambda_i-q^{-2}\lambda_{i+1}} & > & 0, \\
\label{ll-2}
\iota_i\frac{q^2\lambda_{i-1}-\lambda_i}{\lambda_i-\lambda_{i+1}} & > & 0
\end{eqnarray}
for any $i=1,...,n$ and $(\lambda_1,...,\lambda_n)\in {\cal Q}$.
\end{thm}

\begin{pf}
It is clear from \eqref{zl-1}-\eqref{zl-2} that if the conditions
\eqref{ll-1}-\eqref{ll-2} are not satisfied, there is no scalar product
on $W$ that can make it into a unitarizable $\usu$-module.

Suppose that the conditions \eqref{ll-1}-\eqref{ll-2} are satisfied.
By Proposition \ref{Gen-1}, we can identify $W$ with a subspace in
$Func(X)_q$ generated over $Hol(X)_q^{\pm}$ by a function $\delta_{\chi}
\in Func(X)_q^{inv}$ that takes the value $1$ at $\hat{\chi}$ and $0$ at 
any other point in $\cal Q$. The invariant integral \eqref{int-rho}
defines a scalar product on $W$ by $(f,g)=\nu_W(g^*f).$ It is clear
that it defines an irreducible $*$-representation of $Func(X)_q$ and,
hence, of $\usu$, in $W$.
\end{pf}

\begin{exmp}
\label{Gen-5}
Let $U_q{\frak g}_0=U_q{\frak su}(2,1)$, so that $\iota_1=-1$ and $\iota_2=1$.
Then the conditions \eqref{ll-1}-\eqref{ll-2} imply that either
\begin{equation}
\label{21-1}
\lambda_1\geq q^{-1}d_0, \; \lambda_1\geq q^{-2}\lambda_2\geq qc_0
\end{equation}
or
\begin{equation}
\label{21-2}
\lambda_1\leq q^{-1}d_0, \; \lambda_1\leq q^{-2}\lambda_2\leq qc_0.
\end{equation}

It is clear that the lattice $\{(\lambda_1q^k,\lambda_2q^m)\}_{k,m=-\infty}
^{+\infty}$ must be truncated in both horizontal and vertical directions.
It is possible only if either $\zeta_i$ or $\hat{\zeta}_i$ annihilate
a common eigen-vector of $x_1$ and $x_2$. It means that its eigen-value
$(\mu_1,\mu_2)$ must belong to the boundary of the region described by 
either \eqref{21-1} or \eqref{21-2}. Then, if $(\mu_1q^{-2},\mu_2)$ is no
longer in the union of both regions, $\zeta_1$ annihilates the vector,
if $(\mu_1q^2,\mu_2)$ is not in the union of the two regions, 
$\hat{\zeta}_1$ does, and similarly for $\zeta_2$ with $\hat{\zeta}_2$.

It is straightforward to check that only in three cases we get irreducible
$*$-irreducible representations of $U_q{\frak su}(2,1)$. Namely,
\begin{enumerate}
\item when $d_0\leq q^{-2}c_0$, $\lambda_1=q^{-1}d_0$, and $\lambda_2=qc_0$,
we get a highest weight representation, and all the eigen-values of $x_1,x_2$ 
lie in the region described by \eqref{21-2}. It belongs to the degenerate
holomorphic discrete series.
\item when $d_0\geq q^{-2}c_0$, $\lambda_1=q^{-1}d_0$, and $\lambda_2=qc_0$,
we get a lowest weight representation, and all the eigen-values of $x_1,x_2$
lie in the region described by \eqref{21-1}. It belongs to the degenerate
anti-holomorphic discrete series.
\item when $1<\frac{d_0}{c_0}<q^{-4}$, $\lambda_2=qc_0$, and $\lambda_1$ is
such that $(\lambda_1,\lambda_2)$ belongs to the region described by
\eqref{21-1}, while $(\lambda_1q^2,\lambda_2)$ belongs to the region
described by \eqref{21-2}. It belongs to the degenerate complimentary series.
\end{enumerate}

Again, we observe the same 'tunnel effect' that the complimentary series
representations arise in a situation when the set $\cal Q$ of the 
eigen-values of $Func(X)_q^{inv}$ can 'jump' from the set that corresponds
to the holomorphic discrete series representations onto another one that
corresponds to the anti-holomorphic ones. In fact, it is clear that the
regions described by \eqref{21-1} and \eqref{21-2} are nothing but the
projections on $T^*\subset SU(2,1)^*$ of the symplectic leaves of the 
minimal non-zero dimension.
\end{exmp}

\section{Appendix: Holomorphic Realization Of Some Representations Of 
$U_q{\frak su}(1,1)$}

Using the quantum polarizations described above, one can obtain realizations
of the discrete series representations and some of the strange series
representations of $U_q{\frak su}(1,1)$ in spaces of holomorphic
functions. First it was done in \cite{K2}. 

The following propositions are results of straightforward computations.

\begin{prop}
\label{a-1}
\rom{(1)} For any discrete series representation $T_{l,\epsilon}^+$ of
highest weight \rom{(}where $l\leq -\frac{1}{2}$\rom{)}, there exists a 
vector-function $\Theta_+(\lambda)$, holomorphic in the unit disc 
$|\lambda|<1$ and taking values in the space of $T_{l,\epsilon}^+$, such that
\[ \zeta\Theta_+(\lambda)=\lambda\Theta_+(\lambda) \]
and the scalar product of two such functions is equal to
\[ (\Theta_+(\lambda),\Theta_+(\mu))=\frac{1}
{(\lambda\bar{\mu};q^2)_{-2l}}, \]
where $(a;t)_{\alpha}=\frac{(a;t)_{\infty}}{(at^{\alpha};t)_{\infty}}$
and $(a;t)_{\infty}=\prod_{k=0}^{\infty}(1-at^k)$.

\rom{(2)} For any strange series representation $T_{\alpha+\frac{1\pi i}{\log 
q},\epsilon}$, where $2\alpha+1\in {\Bbb N}$, there exists a vector-function 
$\Theta(\lambda)$, holomorphic in the annulus $q^{2\alpha+1}<|\lambda|<1$ 
and taking values in the space of the representation, such that
\[ \zeta\Theta(\lambda)=\lambda\Theta(\lambda) \]
and the scalar product of two such functions is equal to
\[ (\Theta(\lambda),\Theta(\mu))=\frac{(q^{2(\alpha+1-\epsilon)};q^2)
_{\infty}}{(q^{-2(\alpha+\epsilon)};q^2)_{\infty}}\,
_1\Psi_1\left(\begin{array}{c}
q^{-2(\alpha+\epsilon)} \\ q^{2(\alpha+1-\epsilon)}
\end{array};q^2,\lambda\bar{\mu}\right), \]
where
\[ _1\Psi_1\left(\begin{array}{c}
a \\ b
\end{array};t,x\right)=\sum_{k=-\infty}^{+\infty}\frac{(a;t)_k}{(b;t)_k}x^k \]
is the Ramanujan's psi-function.
\end{prop}

Consider the map from the space of $T_{l,\epsilon}^+$ into the space of
functions holomorphic in the unit disc given by
\[ \theta_+:f\mapsto (f,\Theta_+(\bar{\lambda})). \]
Similarly, consider the map from the space of $T_{\alpha+\frac{2\pi i}{\log q},
\epsilon}$ into the space of functions 
holomorphic in the annulus $q^{2\alpha+1}<|\lambda|<1$ given by
\[ \theta:f\mapsto (f,\Theta(\bar{\lambda})). \]

Thus, we get an action of $U_q{\frak su}(1,1)$ on the functions
holomorphic in the respective domain. It is given by $q$-difference
operators. The precise formulas are given in \cite{K2}.

Now we are looking for a measure in the respective domain that establishes
an isomorphism between the corresponding Hilbert spaces -- the space of
the representation and the space of holomorphic functions.

\begin{prop}
\label{a-2}
\rom{(1)} Consider the measure $d\nu_{l,\epsilon}^+$ on the unit disc given by
\[ d\nu_{l,\epsilon}^+=(1-q^{-2(2l-1)})\sum_{k=0}^{\infty}q^{2k}
\frac{(q^{2(k+1)};q^2)_{\infty}}{(q^{2(k-2l-1)};q^2)_{\infty}}
\delta_{|\lambda|=q^k}d\lambda d\bar{\lambda}, \]
when $l<-\frac{1}{2}$, and by
\[ d\nu_{-\frac{1}{2},\epsilon}^+=\delta_{|\lambda|<1}d\lambda 
d\bar{\lambda}, \]
when $l=-\frac{1}{2}$, where $\delta_{|\lambda|=t}$ is the $\delta$-function 
on a circle $|\lambda|=t$. Then $\theta_+$ becomes an isomorphism between the
space of $T_{l,\epsilon}^+$ and the Hilbert space of functions in the unit 
disc with the measure $d\nu_{l,\epsilon}^+$ that are holomorphic in
the interior of the unit disc and continuous in its closure.

\rom{(2)} Consider the measure $d\nu_{\alpha+\frac{2\pi i}
{\log q},\epsilon}$, where $2\alpha+1\in {\Bbb N}$, on the annulus 
$q^{2\alpha+1}\leq |\lambda|\leq 1$ given by
\[ d\nu_{\alpha+\frac{2\pi i}{\log q},\epsilon}=\sum_{k=0}^{2\alpha+1}
q^{2k(\alpha+1-\epsilon)}\frac{(q^{-2(2\alpha+1)};q^2)_k}{(q^2;q^2)_k}
\delta_{|\lambda|=q^k}d\lambda d\bar{\lambda}. \]
Then $\theta$ becomes an isomorphism between the space of $T_{\alpha+
\frac{2\pi i}{\log q},\epsilon}$ and the Hilbert space of the functions
in the above annulus with the measure $d\mu_{\alpha+\frac{2\pi i}{\log q},
\epsilon}$ that are holomorphic in the interior of the annulus and
continuous in it closure.
\end{prop}

Thus, we obtain the realizations of the holomorphic discrete series
representations in the holomorphic functions in the unit disc. And the
realization of some strange series representations in the holomorphic 
functions in an annulus. Similarly, one can obtain a realization of the
anti-holomorphic discrete series representations in the anti-holomorphic
functions in the unit disc (or in the holomorphic functions in the
domain outside the unit disc).

\end{document}